\documentclass[twoside,11pt]{article}

\usepackage{amssymb,amsmath}

\newenvironment{proof}{\begin{trivlist}\item[]{\it
Proof.}}{\hfill$\square$\end{trivlist}}
\newenvironment{proofof}{\noindent{\it Proof of
the theorem.}}{\hfill$\square$\\\mbox{}}

\newenvironment{theorem}{\begin{trivlist}\item[]{\bf\noindent Theorem. }\it}
{\end{trivlist}}

\newtheorem{lemma}{Lemma}
\begin{document}
\title{Simple proof of Chebotar\"ev's theorem on roots of 
unity\footnote{Keywords and phrases: roots of unity, non-zero minors}
\footnote{2000 Mathematics Subject classification: 11T22 (primary),
42A99, 11C20 (secondary)}}
\author{P. E. Frenkel}

\date{}
\newcommand{\Q}{\mathbf Q}
\newcommand{\F}{\mathbf F}
\newcommand{\Z}{\mathbf Z}
\newcommand{\E}{\mathcal E}
\newcommand{\C}{\mathbf C}
\newcommand{\1}{\mathbf 1}
\newcommand{\n}{\mathrm {nilpotent}}

\maketitle

\begin{abstract}  We give  a  simple proof of Chebotar\"ev's theorem: 
Let $p$ be  a prime and $\omega $ a  primitive $p$th root of unity. Then all minors of the matrix 
$\left(\omega^ {ij}\right)_{i,j=0}^{p-1}$ are non-zero.
\end{abstract}

Let $p$ be   a prime   and $\omega  $  a primitive $p$th root of unity.
We write $\F_p$ for the field with $p$ elements.
In 1926, Chebotar\"ev proved the following theorem (see \cite{SL}):

\begin{theorem} For any sets $I,J\subseteq \F_p$ 
with equal cardinality, the matrix
$(\omega^{ij})_{i\in I,j\in J}$ has non-zero determinant. 
\end{theorem}

Several independent proofs have been 
 given, including ones by 
Dieudonn\'e \cite{D}, Evans and Isaacs \cite{EI}, and Terence Tao \cite{T}.  Tao  points out that  the theorem is equivalent to
the inequality $|{\mathrm {supp}} f|+|{\mathrm {supp}} \hat f|\geq p+1$  holding for any function $0\not\equiv f:\F_p\to \C$   
and its Fourier transform $\hat f$, 
a fact also discovered independently by Andr\'as Bir\'o.
Bir\'o posed this as Problem 3 of the 1998 Schweitzer Competition. The proof I
gave in the competition 
(the one in the present article) is published in Hungarian in 
\cite[pp. 53--54.]{matlap}.
It was also discovered (as part of a more general investigation) by Daniel
Goldstein, Robert M. Guralnick and I. M. Isaacs \cite[Section 6]{GGI}.

The proof is based on the following two lemmas. Lemma~\ref{1}
is covered by \cite[Chapter 1]{W}, but we include a proof for the sake of completeness.

\begin{lemma}\label{1}
$\Z[\omega]/(1-\omega)=\F_p.$
\end{lemma}

\begin{proof}
Let $\Omega$ be   an indeterminate and let
$\Phi_p(\Omega)=1+\Omega+\dots+\Omega^{p-1}$ be  the minimal polynomial of 
the algebraic integer $\omega$. Consider the surjective ring homomorphisms
$$\Z[\Omega]\to\Z[\Omega]/(\Phi_p(\Omega))=\Z[\omega], \qquad \Omega\mapsto \omega$$ and
$$\Z[\Omega]\to\Z[\Omega]/(1-\Omega,p)=\F_p, \qquad\Omega\mapsto 1.$$
%
The latter kernel 
contains the former one 
 since $\Phi_p(\Omega)\equiv p \mod (1-\Omega)$.
Therefore, the latter homomorphism 
factors through the former one 
via a surjective homomorphism $\Z[\omega]\to\F_p$ whose kernel is
the ideal  $$(1-\Omega, p)/(\Phi_p(\Omega))=(1-\omega, p)=(1-\omega),$$ the
last equality following from 
$p\equiv\Phi_p(\omega)=0\mod (1-\omega)$.
\end{proof}

\begin{lemma}\label{2}
 Let 
 $0\not\equiv g(x)\in \F_p[x]$ be  a polynomial of degree $<p$. Then
the multiplicity of any element $0\neq a\in \F_p$ as a  root of
$g(x)$ is strictly less than the number of non-zero coefficients of $g(x)$.
\end{lemma}

\begin{proof}
For $g(x)$ constant, the lemma is obviously true.  Assume that
it is true for any 
$g(x)$ of degree $<k$, with some fixed 
$1\leq k<p$  , and take $g(x)$ of degree 
$k$.  If $g(0)=0$, then $g(x)$ has the same number of non-zero coefficients
and the same multiplicity of vanishing at $a$ as $g(x)/x$ does, so the lemma is true for $g(x)$.  If $g(0)\neq 0$, then the number of non-zero coefficients exceeds the corresponding number for
 the derivative $g'(x)$ by 1, and the multiplicity of vanishing at 
$a$ exceeds that of $g'(x)$  by at most 1.  Now $g'(x)\not\equiv 0$ 
since $g(x)$ is of positive degree $k<p$, so the inequality of the 
lemma holds for $g'(x)$
and therefore also for $g(x)$.
\end{proof}

\begin{proofof}
The theorem is equivalent to saying that if  numbers $a_j\in \Q(\omega)$ 
$(j\in J)$
satisfy 
$\sum_{j\in J} a_j\omega^{ij}=0$ for all $i\in I$, then all $a_j$ must be zero.
In fact, we may clearly assume that $a_j\in \Z[\omega]$.
The above equalities mean that the polynomial $$g(x)=\sum_{j\in J} a_j x^j\in \Z[\omega][x]$$
vanishes at $\omega^i$ for all $i\in I$.  So $g(x)$ is divisible by
$\prod_{i\in I}(x-\omega^i)$.  Applying the homomorphism
$\Z[\omega]\to \Z[\omega]/(1-\omega)=\F_p$ to the coefficients of $g(x)$ 
we get a polynomial $\bar g(x)\in \F_p[x]$
that is divisible by $(x-1)^{|I|}$.  On the other hand, 
$\bar g(x)$ has at most $|J|$  non-zero coefficients.   
As $|I|=|J|$, we deduce from Lemma~\ref{2} that 
$\bar g(x)\equiv 0$.  This means that all $a_j$ are divisible by $1-\omega$. 
We may divide all of them by $1-\omega$ and iterate the argument. This leads to
{\it descente infinie} unless all  $a_j$ are zero.
\end{proofof}

\section*{Acknowledgements}
The author is  partially supported by OTKA grants T 042769 and T 046365.
The author wishes to thank Andr\'as Bir\'o for calling his attention to the 
subject and  for useful discussions.

\bigskip
\noindent
{\bf Address.} Mathematics Institute, 
Budapest University of Technology and Economics, Egry J. u. 1., Budapest, 
1111 Hungary.  
E-mail:  frenkelp@renyi.hu

\end{document}